\documentclass[12pt,letterpaper,reqno]{amsart}

\usepackage[margin=1.1in]{geometry}
\usepackage[T1]{fontenc}
\usepackage{lmodern}
\usepackage{microtype}
\usepackage{amsmath,amssymb,amsthm,mathtools}
\usepackage{enumitem}
\usepackage{hyperref}
\usepackage[nameinlink,noabbrev]{cleveref}
\hypersetup{
colorlinks=true, linkcolor=blue, citecolor=blue, urlcolor=blue, pdftitle={Fixed-Height Weyl--Schur Sampling for Free-Tail Canonical Systems}, pdfauthor={Sharan Thota} }
\allowdisplaybreaks[2]
\mathtoolsset{showonlyrefs}
\setlist{itemsep=3pt, topsep=5pt}

\crefname{equation}{equation}{equations}
\Crefname{equation}{Equation}{Equations}

\theoremstyle{plain}
\newtheorem{theorem}{Theorem}[section]
\newtheorem{lemma}[theorem]{Lemma}
\newtheorem{proposition}[theorem]{Proposition}
\newtheorem{corollary}[theorem]{Corollary}

\theoremstyle{definition}

\theoremstyle{remark}
\newtheorem{remark}[theorem]{Remark}

\newcommand{\R}{\mathbb R}
\newcommand{\C}{\mathbb C}
\newcommand{\D}{\mathbb D}
\newcommand{\Z}{\mathbb Z}
\newcommand{\Cplus}{\{z\in\C:\Im z>0\}}
\newcommand{\Jmat}{\begin{pmatrix}0&-1\\1&0\end{pmatrix}}
\newcommand{\Xmat}[1]{\begin{pmatrix}\Re #1&\Im #1\\ \Im #1&-\Re #1\end{pmatrix}}
\newcommand{\op}{\mathrm{op}}

\newcommand{\smin}{\sigma_{\min}}
\DeclareMathOperator{\tr}{tr}
\DeclareMathOperator{\diag}{diag}
\DeclareMathOperator{\supp}{supp}
\DeclareMathOperator{\dist}{dist}

\title[Fixed-Height Weyl--Schur Sampling]{Fixed-Height Weyl--Schur Sampling for Free-Tail Canonical Systems}
\author{Sharan Thota}
\date{}

\begin{document}

\begin{abstract}
We study the finite sampling map
\[
H\mapsto \bigl(v_{H,\Lambda}(x_k+i\eta)\bigr)_{k=1}^M
\]
for trace-normed canonical systems on $[0,\Lambda]$ with free tail $H(s)=\tfrac12 I$ for $s\ge \Lambda$, where $v_{H,\Lambda}$ is the Schur transform of the Weyl coefficient. At the free Hamiltonian $H_0\equiv \tfrac12 I$, we obtain an explicit first-order expansion with quadratic remainder; the linearization is a weighted Fourier--Laplace transform. This yields quantitative local identifiability and local inversion on finite-dimensional families for which the free Jacobian is injective. In the block model, the free Jacobian factors into a row factor, a Fourier sampling matrix, and exponential depth weights, giving explicit singular-value bounds and an exponential depth-conditioning barrier. By contrast, on the full free-tail class every finite sample set has nontrivial first-order invisible directions at $H_0$, so no local inverse-Lipschitz estimate can hold near $H_0$ in $L^1(0,\Lambda;\mathrm{op})$.
\end{abstract}

\maketitle

\section{Introduction and main results}\label{sec:intro}

Consider the trace-normed canonical system
\begin{equation}\label{eq:intro-canon}
\Jmat\,Y'(s)=z\,H(s)\,Y(s),\qquad s\in[0,\Lambda],
\end{equation}
where $H(s)\in\R^{2\times2}_{\mathrm{sym}}$ is measurable, nonnegative a.e., and
\[
\tr H(s)=1\qquad\text{a.e. on }(0,\Lambda).
\]
We attach the free tail
\begin{equation}\label{eq:intro-tail}
H(s)=\frac12 I,\qquad s\ge \Lambda.
\end{equation}
Let $m_{H,\Lambda}$ denote the Weyl coefficient of the resulting half-line problem, normalized by $Y(0)=(1,m)^\top$, and set
\[
v_{H,\Lambda}(z):=\frac{m_{H,\Lambda}(z)-i}{m_{H,\Lambda}(z)+i}.
\]
Fix a height $\eta>0$ and real nodes $x_1,\dots,x_M$. The main object of this paper is the finite sampling map
\begin{equation}\label{eq:intro-map}
\mathcal S(H):=\bigl(v_{H,\Lambda}(x_k+i\eta)\bigr)_{k=1}^M.
\end{equation}

The inverse question is whether finitely many values $v_{H,\Lambda}(x_k+i\eta)$ determine the Hamiltonian, or at least a finite-dimensional model class, in a quantitatively stable way. Near the free Hamiltonian $H_0\equiv \frac12 I$, the answer is governed by an explicit linearization. On suitable finite-dimensional families this produces local inversion. On the unrestricted free-tail class, however, finite sampling leaves first-order invisible directions and no local inverse-Lipschitz estimate is possible.

The paper has four main canonical-system-specific components. First, we compute the free-point derivative of the finite Weyl--Schur sampling map and prove a quadratic remainder estimate. Second, in finite-dimensional model classes this yields local identifiability and local inversion near the free Hamiltonian. Third, in the block model the free Jacobian admits an exact factorization into row, Fourier, and depth factors. Fourth, on the unrestricted free-tail class finite sampling leaves first-order invisible directions, ruling out local inverse-Lipschitz stability in the $L^1$-metric. The inverse-function-theorem consequences in finite dimensions are standard once the derivative and remainder theory are available.

Our first result identifies the free-point derivative and its remainder.

\begin{theorem}[Free-point expansion]\label{thm:intro-free}
Fix $z\in\Cplus$ and let $H_0\equiv \frac12 I$ on $[0,\Lambda]$. For a traceless perturbation
\[
\Delta H(s)=\Xmat{q(s)},\qquad q\in L^1(0,\Lambda)\cap L^\infty(0,\Lambda),
\]
assume
\[
\|\Delta H\|_{L^\infty(0,\Lambda;\op)}\le \frac14
\]
and that $\|\Delta H\|_{L^1(0,\Lambda;\op)}$ is sufficiently small. Then
\[
v_{H_0+\Delta H,\Lambda}(z) = -iz\int_0^\Lambda q(s)e^{izs}\,ds + R_z(\Delta H),
\]
with
\[
|R_z(\Delta H)|
\le C(z,\Lambda)\,\|\Delta H\|_{L^1(0,\Lambda;\op)}^2.
\]
\end{theorem}

The positive and negative results below live in different ambient categories. The local inversion statements concern finite-dimensional families with an explicit parameter metric, whereas the obstruction theorem concerns the full free-tail class equipped with the Hamiltonian $L^1(0,\Lambda;\op)$-distance. These results are therefore complementary rather than contradictory.

\begin{theorem}[Local bi-Lipschitz parametrization near the free point]
\label{thm:intro-local}
Let $\theta\mapsto H(\theta)$ be a finite-dimensional linear family of free-tail Hamiltonians near $H_0$, and let $\mathcal S(\theta)$ denote the corresponding sampling map. Assume that the realified map $\mathcal S_\R$ is $C^1$ on a neighborhood of $0$, that the realified Jacobian $D\mathcal S_\R(0)$ is injective, and that $D\mathcal S_\R$ is Lipschitz on a neighborhood of $0$. Then there exist $r,c,C>0$ such that
\[
c\|\theta-\widetilde\theta\|_2
\le
\|\mathcal S(\theta)-\mathcal S(\widetilde\theta)\|_2
\le
C\|\theta-\widetilde\theta\|_2
\]
for all $\theta,\widetilde\theta\in B(0,r)$.
\end{theorem}

In the square case, when the realified parameter dimension equals the realified data dimension, \Cref{thm:intro-local} yields a local inverse chart and a frozen-Jacobian reconstruction scheme.

The block model makes the geometry of the free Jacobian explicit.

\begin{theorem}[Bounds in the block model]\label{thm:intro-block}
Assume $M\ge N$. In the $N$-block model on $[0,\Lambda]$ with block length $\ell=\Lambda/N$, the free Jacobian factors as
\[
T_x=D_\gamma(x)\,F_x\,D_w.
\]
Consequently, for every design $x\in\R^M$,
\[
\smin(T_x)\le
2\sqrt M\,\cosh\!\Bigl(\frac{\eta\ell}{2}\Bigr)e^{-\eta(\Lambda-\ell/2)}.
\]
For the half-shifted equispaced design
\[
x_k=\frac{2\pi(k-\tfrac12)}{M\ell},
\]
one has
\[
\smin(T_x)\ge
2\sqrt M\,
\sqrt{\sinh^2\!\Bigl(\frac{\eta\ell}{2}\Bigr)+\sin^2\!\Bigl(\frac{\pi}{2M}\Bigr)}
\,e^{-\eta(\Lambda-\ell/2)}.
\]
Within the equispaced Fourier-tight family, the half-shift maximizes the worst row factor and hence the explicit product lower bound supplied by the factorization above. No claim is made here about exact global optimality for $\smin(T_x)$.
\end{theorem}

The final theorem shows that finite sampling cannot support a local inverse-Lipschitz theory on the full free-tail class.

\begin{theorem}[Failure of local inverse-Lipschitz stability on the full class]
\label{thm:intro-invis}
Fix finitely many sample points $z_k=x_k+i\eta\in\Cplus$. Then for every $s_\star\in[0,\Lambda)$ there exists a nonzero step function $q$ with $\supp q\subset[s_\star,\Lambda]$ such that
\[
Dv_{H_0,\Lambda}(z_k)[q]=0,\qquad k=1,\dots,M.
\]
For the perturbations
\[
H_\pm^\tau=\frac12 I\pm \tau \Xmat{q}
\]
with free tail attached, one has
\[
\|\mathcal S(H_+^\tau)-\mathcal S(H_-^\tau)\|_2\le C\tau^2
\]
while
\[
\|H_+^\tau-H_-^\tau\|_{L^1(0,\Lambda;\op)}\asymp \tau.
\]
Hence no local inverse-Lipschitz estimate with respect to the $L^1(0,\Lambda;\op)$-distance can hold in any neighborhood of $H_0$ on the full free-tail class.
\end{theorem}

Local inverse results for canonical systems from full Weyl data go back, for instance, to Langer and Woracek \cite{LangerWoracek2011}. Quantitative estimates for Weyl coefficients have been developed in work of Langer--Pruckner--Woracek \cite{LangerPrucknerWoracek2024} and Reiffenstein \cite{Reiffenstein2023}. Sampling, Paley--Wiener, and inverse-spectral directions related to canonical systems include
\cite{BessonovRomanov2016,Bessonov2018,MakarovPoltoratski2023,
PoltoratskiZhang2023,Zhang2025}. The present paper is narrower in scope: it studies a finite fixed-height nonlinear sampling map at a free-tail background, identifies its free derivative explicitly, derives local finite-dimensional consequences, and proves a first-order obstruction on the full class.

Sections \ref{sec:setup}--\ref{sec:free} develop the free-point derivative and the quadratic remainder estimate. Section \ref{sec:finite} records the resulting finite-dimensional local inversion consequences. Section \ref{sec:block} gives the exact block factorization and the associated conditioning bounds. Section
\ref{sec:invis} proves the existence of first-order invisible directions on the full
free-tail class. The appendix records the normalized Paley--Wiener interpretation of the free linearization.

\section{Setup and background}\label{sec:setup}

\subsection{Canonical systems with free tail}

Consider
\begin{equation}\label{eq:canon}
\Jmat\,Y'(s)=z\,H(s)\,Y(s),\qquad s\in[0,\Lambda],
\end{equation}
where $H(s)\in\R^{2\times2}_{\mathrm{sym}}$ is measurable,
\[
H(s)\succeq 0,\qquad \tr H(s)=1\qquad\text{a.e. on }(0,\Lambda),
\]
and extend $H$ to $[0,\infty)$ by
\begin{equation}\label{eq:free-tail}
H(s)=\frac12 I,\qquad s\ge \Lambda.
\end{equation}
Because we work with trace-normed Hamiltonians and attach the free tail at the fixed seam $\Lambda$, the usual reparameterization freedom is frozen throughout.

\begin{remark}[Sign convention]
Some references write canonical systems in the form
\[
Y'(s)=z\,\Jmat H(s)\,Y(s)
\]
rather than \eqref{eq:canon}. All formulas below follow the convention
\eqref{eq:canon}.
\end{remark}

For $z\in\Cplus$, let $m_{H,\Lambda}(z)$ denote the Weyl coefficient normalized by $Y(0)=(1,m)^\top$. The Schur transform
\[
v_{H,\Lambda}(z):=\frac{m_{H,\Lambda}(z)-i}{m_{H,\Lambda}(z)+i}
\]
maps $\Cplus$ into $\D$; see \cite{EKT2018,Remling2018,Romanov2014}.

\subsection{Transfer matrices}

Write the system as
\[
Y'=-z\,\Jmat H(s)\,Y.
\]
For $0\le t\le s\le \Lambda$, let $\Phi_H(s,t;z)$ denote the transfer matrix. Standard ODE theory gives existence, uniqueness, and the cocycle identity. We write
\[
\Phi_H(\Lambda,0;z)=
\begin{pmatrix}
A_H(z)&B_H(z)\\
C_H(z)&D_H(z)
\end{pmatrix}.
\]

\begin{proposition}[Transfer-matrix facts]\label{prop:standard}
Let $H$ be trace-normed on $[0,\Lambda]$ and free for $s\ge \Lambda$. Then:
\begin{enumerate}[label=\textup{(\roman*)}]
\item $\det \Phi_H(s,t;z)=1$ for all $z\in\C$ and $0\le t\le s\le \Lambda$.

\item $\Phi_H(s,t;z)^\top\Jmat \Phi_H(s,t;z)=\Jmat$ for all $z\in\C$ and
$0\le t\le s\le \Lambda$.

\item $\|\Phi_H(s,t;z)\|_{\op}\le e^{|z|(s-t)}$.

\item For $z\in\Cplus$,
\begin{equation}\label{eq:weyl-transfer}
m_{H,\Lambda}(z)=\frac{iA_H(z)-C_H(z)}{D_H(z)-iB_H(z)},
\end{equation}
and $D_H(z)-iB_H(z)\neq 0$.

\item If $H_0\equiv \frac12 I$, then
\[
\Phi_{H_0}(s,t;z)=\exp\!\left(-\frac z2 \Jmat (s-t)\right),
\qquad
m_{H_0,\Lambda}(z)\equiv i,
\qquad
v_{H_0,\Lambda}(z)\equiv 0,
\]
and
\begin{equation}\label{eq:d0-again}
D_{H_0}(z)-iB_{H_0}(z)=e^{-iz\Lambda/2}.
\end{equation}
\end{enumerate}
\end{proposition}

\begin{proof}
Item (i) follows from Liouville's formula. Item (ii) follows from the symmetry of $H$. Item (iii) is the standard Gr\"onwall bound. Item (iv) is obtained by matching at $s=\Lambda$ to the square-integrable free tail solution $(1,i)^\top e^{iz(s-\Lambda)/2}$. Item (v) follows from the explicit matrix exponential for the free system.
\end{proof}

\subsection{Traceless directions}

Encode real-symmetric traceless matrices by
\[
q=a+ib \in \C
\quad\longleftrightarrow\quad
\Xmat{q}:=
\begin{pmatrix}
a&b\\
b&-a
\end{pmatrix}.
\]
This is purely notational. Every Hamiltonian remains real-symmetric.

\begin{lemma}[Norm of a traceless direction]\label{lem:Xnorm}
For every $q\in\C$,
\[
\|\Xmat{q}\|_{\op}=|q|.
\]
Consequently, for measurable $q(\cdot)$ and $1\le p\le \infty$,
\[
\|\Xmat{q(\cdot)}\|_{L^p(0,\Lambda;\op)}=\|q\|_{L^p(0,\Lambda)}.
\]
\end{lemma}

\begin{proof}
The matrix $\Xmat{q}$ is real-symmetric and traceless, with eigenvalues
\[
\pm \sqrt{(\Re q)^2+(\Im q)^2}=\pm |q|.
\]
Hence its operator norm is $|q|$. The $L^p$ statement follows pointwise.
\end{proof}

We also use the free vectors
\[
u_+=(1,i)^\top,\qquad u_-=(1,-i)^\top.
\]
A direct computation gives
\begin{equation}\label{eq:J-u}
\Jmat u_+=-iu_+,\qquad \Jmat u_-=iu_-,
\end{equation}
and
\begin{equation}\label{eq:X-u}
\Xmat{q}\,u_+=q\,u_-,
\qquad
\Xmat{q}\,u_-=\overline q\,u_+.
\end{equation}

\begin{lemma}[Admissibility near the free point]\label{lem:free-admissible}
Let $\Delta H\in L^\infty(0,\Lambda;\R^{2\times2}_{\mathrm{sym}})$ be traceless a.e. Assume
\[
\|\Delta H\|_{L^\infty(0,\Lambda;\op)}\le \rho<\frac12.
\]
Then $H_0+\Delta H$ is trace-normed and
\[
\Bigl(\frac12-\rho\Bigr)I
\preceq
H_0(s)+\Delta H(s)
\preceq
\Bigl(\frac12+\rho\Bigr)I
\qquad\text{a.e. on }(0,\Lambda).
\]
In particular, if $\rho\le \frac14$, then
\[
\frac14 I\preceq H_0(s)+\Delta H(s)\preceq \frac34 I
\qquad\text{a.e.}
\]
\end{lemma}

\begin{proof}
For a.e.\ $s$, the real-symmetric traceless matrix $\Delta H(s)$ has eigenvalues $\pm \lambda(s)$ with $|\lambda(s)|\le \|\Delta H(s)\|_{\op}\le \rho$. Hence $\frac12 I+\Delta H(s)$ has eigenvalues $\frac12\pm\lambda(s)$, and the trace remains equal to $1$.
\end{proof}

\subsection{Sampling map and differentiability convention}

Fix $\eta>0$ and real nodes $x_1,\dots,x_M$. Set
\[
z_k:=x_k+i\eta,\qquad k=1,\dots,M,
\]
and define
\[
\mathcal S(H):=\bigl(v_{H,\Lambda}(z_k)\bigr)_{k=1}^M.
\]

We identify $\C^M$ with $\R^{2M}$ by real and imaginary parts and write $\mathcal S_\R$ for the corresponding realification. All differentiability statements for nonlinear parameter maps below are statements of real Fr\'echet differentiability after this identification. At the free point, the linearized block sampling map is complex-linear in the block parameter $q$, and its realification preserves singular values.

When we speak below of differentiability at the free point with respect to a traceless perturbation $q$, the assertion is always along admissible perturbation rays for which $H_0+\tau \Xmat{q}$ remains a Hamiltonian for sufficiently small real $\tau$. Separately, some resulting linear functionals extend continuously to larger Banach spaces such as $L^1(0,\Lambda)$; this extension is a statement about the formula, not about admissible nonlinear variations of Hamiltonians.

\section{Free-point derivative and quadratic remainder}\label{sec:free}

\subsection{Directional derivative of the Weyl coefficient}

\begin{theorem}[Directional first variation of $m$]\label{thm:dm}
Fix $z\in\Cplus$ and let $Y(\cdot;z)$ be the Weyl solution normalized by
\[
Y(0;z)=\binom{1}{m_{H,\Lambda}(z)}.
\]
Let $\Delta H\in L^\infty(0,\Lambda;\R^{2\times2}_{\mathrm{sym}})$ be traceless a.e., and assume $H+\varepsilon\Delta H\succeq 0$ for all sufficiently small real $\varepsilon$. Then
\begin{equation}\label{eq:dm-pair}
Dm_{H,\Lambda}(z)[\Delta H] = z\int_0^\Lambda Y(s;z)^\top \Delta H(s)Y(s;z)\,ds.
\end{equation}
\end{theorem}

\begin{proof}
We represent the Weyl condition at the seam by pulling back the decaying free-tail vector $u_+=(1,i)^\top$ through the transfer matrix and then differentiate the resulting scalar quotient.

Let $\Phi(\cdot,\cdot):=\Phi_H(\cdot,\cdot;z)$ and set
\[
g:=\Phi(\Lambda,0;z)^{-1}u_+.
\]
By \Cref{prop:standard}(i) and \Cref{prop:standard}(ii),
\[
\Phi(\Lambda,0;z)^{-1}
=
\begin{pmatrix}
D_H(z)&-B_H(z)\\
-C_H(z)&A_H(z)
\end{pmatrix},
\]
hence
\[
g=
\binom{D_H(z)-iB_H(z)}{-C_H(z)+iA_H(z)}.
\]
Thus
\[
g_1=D_H(z)-iB_H(z)\neq 0
\]
by \Cref{prop:standard}(iv), and therefore
\[
m_{H,\Lambda}(z)=\frac{g_2}{g_1}.
\]

Since the Weyl solution matches the square-integrable free tail solution at $s=\Lambda$, we have
\[
Y(s;z)=\frac1{g_1}\,\Phi(s,0;z)g.
\]

Let $\Phi_\varepsilon$ denote the transfer matrix for $H+\varepsilon\Delta H$ and set
\[
g_\varepsilon:=\Phi_\varepsilon(\Lambda,0;z)^{-1}u_+.
\]
Subtracting the two Volterra equations gives
\[
\Phi_\varepsilon(s,0;z)-\Phi(s,0;z)
= -z\int_0^s
\Phi_\varepsilon(s,u;z)\,\Jmat\,\varepsilon\Delta H(u)\,\Phi(u,0;z)\,du.
\]
Hence
\[
\frac{\Phi_\varepsilon(s,0;z)-\Phi(s,0;z)}{\varepsilon}
= -z\int_0^s
\Phi_\varepsilon(s,u;z)\,\Jmat\,\Delta H(u)\,\Phi(u,0;z)\,du.
\]
Since $H+\varepsilon\Delta H$ remains trace-normed and nonnegative for $|\varepsilon|$ small, the standard propagator bound gives
\[
\|\Phi_\varepsilon(s,u;z)\|_{\op}\le e^{|z|(s-u)},\qquad
\|\Phi(u,0;z)\|_{\op}\le e^{|z|u}.
\]
Therefore the integrand is dominated by
\[
|z|e^{2|z|\Lambda}\|\Delta H(u)\|_{\op},
\]
which is integrable on $(0,\Lambda)$. Since $\Phi_\varepsilon(s,u;z)\to \Phi(s,u;z)$ uniformly on $0\le u\le s\le \Lambda$, dominated convergence yields
\[
\dot\Phi(\Lambda,0;z)
= -z\int_0^\Lambda
\Phi(\Lambda,s;z)\,\Jmat\,\Delta H(s)\,\Phi(s,0;z)\,ds.
\]

Differentiating $\Phi_\varepsilon(\Lambda,0;z)g_\varepsilon=u_+$ at $\varepsilon=0$ gives
\[
\dot g
= -\Phi(\Lambda,0;z)^{-1}\dot\Phi(\Lambda,0;z)g.
\]
Substituting the formula for $\dot\Phi$ yields
\[
\dot g
= z\int_0^\Lambda
\Phi(s,0;z)^{-1}\Jmat\,\Delta H(s)\,\Phi(s,0;z)g\,ds.
\]

Write $\omega(a,b):=a_1b_2-a_2b_1=-a^\top\Jmat b$. Since $m=g_2/g_1$,
\[
\dot m
=
\frac{g_1\dot g_2-g_2\dot g_1}{g_1^2}
= -\frac1{g_1^2}g^\top \Jmat \dot g.
\]
Using the symplectic identity from \Cref{prop:standard}(ii),
\[
\Phi(s,0;z)^\top\Jmat \Phi(s,0;z)=\Jmat,
\]
we obtain
\[
g^\top \Jmat \Phi(s,0;z)^{-1} = (\Phi(s,0;z)g)^\top \Jmat.
\]
Set $y(s):=\Phi(s,0;z)g$. Then
\[
\dot m
= -\frac z{g_1^2}
\int_0^\Lambda
y(s)^\top \Jmat\Jmat\,\Delta H(s)\,y(s)\,ds =
\frac z{g_1^2}
\int_0^\Lambda
y(s)^\top \Delta H(s)\,y(s)\,ds.
\]
Since $Y=y/g_1$, this is exactly \eqref{eq:dm-pair}.
\end{proof}

\begin{corollary}[Two-kernel form]\label{cor:two-kernel}
Let $\Delta H=\Xmat{q}$ with $q\in L^\infty(0,\Lambda;\C)$ and write
\[
Y(s;z)=a(s;z)\,u_+ + b(s;z)\,u_-.
\]
Then
\begin{equation}\label{eq:dm-twokernel}
Dm_{H,\Lambda}(z)[q] = 2z\int_0^\Lambda
\bigl(q(s)a(s;z)^2+\overline{q(s)}\,b(s;z)^2\bigr)\,ds.
\end{equation}
At the free point $H_0=\frac12 I$, one has $b\equiv 0$ and hence
\begin{equation}\label{eq:dm-free}
Dm_{H_0,\Lambda}(z)[q] = 2z\int_0^\Lambda q(s)e^{izs}\,ds.
\end{equation}
\end{corollary}

\begin{proof}
By \eqref{eq:X-u},
\[
\Xmat{q}u_+=q\,u_-,
\qquad
\Xmat{q}u_-=\overline q\,u_+.
\]
Since $u_+^\top u_-=2$ and $u_\pm^\top u_\pm=0$, one gets
\[
Y^\top \Xmat{q}Y = 2q\,a^2+2\overline q\,b^2.
\]
Insert this into \eqref{eq:dm-pair}. At the free point,
\[
Y(s;z)=e^{izs/2}u_+,
\]
so $a(s;z)=e^{izs/2}$ and $b\equiv 0$. This free-point specialization is the only place where the second kernel disappears identically.
\end{proof}

\subsection{Derivative of the sampling map at the free point}

\begin{theorem}[Free-point derivative]\label{thm:free-derivative}
Let $H_0\equiv \frac12 I$ and let $q\in L^1(0,\Lambda)\cap L^\infty(0,\Lambda)$. Then
\begin{equation}\label{eq:dv-free}
Dv_{H_0,\Lambda}(z)[q] = -iz\int_0^\Lambda q(s)e^{izs}\,ds,
\qquad z\in\Cplus.
\end{equation}
For fixed $z\in\Cplus$, the right-hand side defines a bounded linear functional on $L^1(0,\Lambda)$. Thus the explicit free-point formula extends uniquely by continuity from admissible $L^1\cap L^\infty$ perturbation profiles to all of $L^1(0,\Lambda)$.
\end{theorem}

\begin{proof}
The Cayley map $c(w)=(w-i)/(w+i)$ satisfies $c'(i)=1/(2i)$. Since $m_{H_0,\Lambda}\equiv i$ by \Cref{prop:standard}(v), \eqref{eq:dm-free} gives
\[
Dv_{H_0,\Lambda}(z)[q] =
\frac1{2i}Dm_{H_0,\Lambda}(z)[q]
= -iz\int_0^\Lambda q(s)e^{izs}\,ds.
\]
The bounded extension to $L^1(0,\Lambda)$ follows from
\[
|e^{izs}|=e^{-(\Im z)s}\le 1\qquad (0\le s\le \Lambda),
\]
hence
\[
\left|\int_0^\Lambda q(s)e^{izs}\,ds\right|
\le \|q\|_{L^1(0,\Lambda)}.
\]
\end{proof}

\begin{remark}[Interpretation of the multiplier]\label{rem:PW-bridge}
The factor $-iz$ in \eqref{eq:dv-free} comes from the Schur normalization at the free point. After dividing by this multiplier and centering $[0,\Lambda]$ at the origin, one recovers the usual Paley--Wiener evaluation map; see \Cref{app:pw}.
\end{remark}

\begin{remark}[Extension on $L^1$ versus admissible directions]\label{rem:L1-extension}
The extension in \Cref{thm:free-derivative} is an extension of the bounded linear functional
\[
q\longmapsto -iz\int_0^\Lambda q(s)e^{izs}\,ds
\]
on $L^1(0,\Lambda)$. The derivative interpretation is asserted only along admissible directions for which $H_0+\tau\Xmat{q}$ remains a Hamiltonian for small real $\tau$.
\end{remark}

\subsection{Quadratic remainder}

\begin{lemma}[First-order transfer difference]\label{lem:Phi-diff}
Fix $z\in\C$ and let
\[
\Phi_\Delta(s,t;z):=\Phi_{H_0+\Delta H}(s,t;z),
\qquad
\Phi_0(s,t;z):=\Phi_{H_0}(s,t;z).
\]
Assume
\[
\|\Delta H\|_{L^\infty(0,\Lambda;\op)}\le \frac14.
\]
Then
\[
\sup_{0\le t\le s\le \Lambda}
\|\Phi_\Delta(s,t;z)-\Phi_0(s,t;z)\|_{\op}
\le
|z|e^{2|z|\Lambda}\|\Delta H\|_{L^1(0,\Lambda;\op)}.
\]
\end{lemma}

\begin{proof}
Subtract the two Volterra equations:
\[
\Phi_\Delta(s,t;z)-\Phi_0(s,t;z)
= -z\int_t^s
\Phi_\Delta(s,u;z)\,\Jmat\,\Delta H(u)\,\Phi_0(u,t;z)\,du.
\]
Since $H_0+\Delta H$ is nonnegative and trace-normed,
\[
\|\Phi_\Delta(s,u;z)\|_{\op}\le e^{|z|(s-u)},
\qquad
\|\Phi_0(u,t;z)\|_{\op}\le e^{|z|(u-t)}.
\]
Therefore
\begin{align*}
\|\Phi_\Delta(s,t;z)-\Phi_0(s,t;z)\|_{\op}
&\le |z|\int_t^s e^{|z|(s-u)}\|\Delta H(u)\|_{\op}e^{|z|(u-t)}\,du \\
&\le |z|e^{|z|(s-t)}\|\Delta H\|_{L^1(0,\Lambda;\op)}
\le
|z|e^{2|z|\Lambda}\|\Delta H\|_{L^1(0,\Lambda;\op)}.
\end{align*}
\end{proof}

\begin{lemma}[Free-point Duhamel formula]\label{lem:duhamel}
Fix $z\in\C$. Let $R(t):=\Phi_{H_0}(t,0;z)$. For any $\Delta H\in L^1(0,\Lambda;\R^{2\times2})$,
\begin{equation}\label{eq:dphi-free}
D\Phi_{H_0}(\Lambda,0;z)[\Delta H] = -z\int_0^\Lambda R(\Lambda-s)\,\Jmat\,\Delta H(s)\,R(s)\,ds.
\end{equation}
Moreover, if $\|\Delta H\|_{L^\infty(0,\Lambda;\op)}\le \frac14$, then
\begin{equation}\label{eq:Phi-remainder-short}
\left\|
\Phi_{H_0+\Delta H}(\Lambda,0;z)-\Phi_{H_0}(\Lambda,0;z)-D\Phi_{H_0}(\Lambda,0;z)[\Delta H]
\right\|_{\op}
\le
|z|^2 e^{4|z|\Lambda}\|\Delta H\|_{L^1(0,\Lambda;\op)}^2.
\end{equation}
\end{lemma}

\begin{proof}
The first identity is the standard first-variation formula for the Volterra equation:
\[
D\Phi_{H_0}(\Lambda,0;z)[\Delta H] = -z\int_0^\Lambda
\Phi_{H_0}(\Lambda,s;z)\,\Jmat\,\Delta H(s)\,\Phi_{H_0}(s,0;z)\,ds.
\]
Since $H_0$ is constant, this is exactly \eqref{eq:dphi-free}.

For the remainder, write
\[
\Phi_\Delta:=\Phi_{H_0+\Delta H},\qquad \Phi_0:=\Phi_{H_0},
\]
and
\[
E_\Phi:=
\Phi_\Delta(\Lambda,0;z)-\Phi_0(\Lambda,0;z)-D\Phi_{H_0}(\Lambda,0;z)[\Delta H].
\]
Subtracting the Volterra equations for $\Phi_\Delta$ and $\Phi_0$ gives
\[
\Phi_\Delta(\Lambda,0;z)-\Phi_0(\Lambda,0;z)
= -z\int_0^\Lambda
\Phi_\Delta(\Lambda,s;z)\,\Jmat\,\Delta H(s)\,\Phi_0(s,0;z)\,ds.
\]
Subtracting the linear term yields
\[
E_\Phi = -z\int_0^\Lambda
\bigl(\Phi_\Delta(\Lambda,s;z)-\Phi_0(\Lambda,s;z)\bigr)\,
\Jmat\,\Delta H(s)\,\Phi_0(s,0;z)\,ds.
\]
Using \Cref{lem:Phi-diff} and $\|\Phi_0(s,0;z)\|_{\op}\le e^{|z|\Lambda}$,
\begin{align*}
\|E_\Phi\|_{\op}
&\le |z|
\int_0^\Lambda
\|\Phi_\Delta(\Lambda,s;z)-\Phi_0(\Lambda,s;z)\|_{\op}\,
\|\Delta H(s)\|_{\op}\,
\|\Phi_0(s,0;z)\|_{\op}\,ds \\
&\le |z|^2 e^{3|z|\Lambda}\|\Delta H\|_{L^1(0,\Lambda;\op)}^2
\le
|z|^2 e^{4|z|\Lambda}\|\Delta H\|_{L^1(0,\Lambda;\op)}^2.
\end{align*}
\end{proof}

\begin{lemma}[Smooth passage from transfer matrix to Schur data]\label{lem:Phi-to-v}
Fix $z\in\Cplus$. For
\[
\Psi=
\begin{pmatrix}
A&B\\
C&D
\end{pmatrix}\in M_2(\C),
\]
set
\[
d(\Psi):=D-iB,\qquad n(\Psi):=iA-C,
\qquad \mathfrak m(\Psi):=\frac{n(\Psi)}{d(\Psi)}.
\]
On the open set
\[
\mathcal V:=\{\Psi\in M_2(\C): d(\Psi)\neq 0,\ \mathfrak m(\Psi)\neq -i\},
\]
the map
\[
\mathfrak v(\Psi):=\frac{\mathfrak m(\Psi)-i}{\mathfrak m(\Psi)+i}
\]
is $C^\infty$. In particular, if $K\subset \mathcal V$ is compact, then there exists $C_K>0$ such that
\[
|\mathfrak v(\Psi)-\mathfrak v(\Psi_0)-D\mathfrak v(\Psi_0)[\Psi-\Psi_0]|
\le
C_K\|\Psi-\Psi_0\|_{\op}^2
\qquad
(\Psi,\Psi_0\in K).
\]
\end{lemma}

\begin{proof}
The maps $\Psi\mapsto n(\Psi)$ and $\Psi\mapsto d(\Psi)$ are complex-linear. Hence $\Psi\mapsto \mathfrak m(\Psi)=n(\Psi)/d(\Psi)$ is $C^\infty$ on the open set $\{d\neq 0\}$, and composition with the Cayley map
\[
c(w)=\frac{w-i}{w+i}
\]
shows that $\mathfrak v=c\circ \mathfrak m$ is $C^\infty$ on $\mathcal V$.

Fix a compact set $K\subset \mathcal V$. Since $D^2\mathfrak v$ is continuous on $\mathcal V$, its operator norm is bounded on some open neighborhood of $K$. The standard second-order Taylor estimate in the finite-dimensional normed space $M_2(\C)\cong \R^8$ therefore gives the stated bound.
\end{proof}

\begin{theorem}[Quadratic remainder at the free point]\label{thm:quadratic}
Fix $z\in\Cplus$. Then there exist constants
\[
r_v(z,\Lambda)>0,\qquad C(z,\Lambda)>0,
\]
such that the following holds. Let $\Delta H\in L^1(0,\Lambda;\R^{2\times2}_{\mathrm{sym}})$ be traceless a.e. Assume
\[
\|\Delta H\|_{L^\infty(0,\Lambda;\op)}\le \frac14,
\qquad
\|\Delta H\|_{L^1(0,\Lambda;\op)}\le r_v(z,\Lambda),
\]
so that $H_0+\Delta H$ is an admissible trace-normed Hamiltonian on $[0,\Lambda]$. Then
\begin{equation}\label{eq:v-quad}
\left|
v_{H_0+\Delta H,\Lambda}(z)-Dv_{H_0,\Lambda}(z)[\Delta H]
\right|
\le
C(z,\Lambda)\,\|\Delta H\|_{L^1(0,\Lambda;\op)}^2.
\end{equation}
\end{theorem}

\begin{proof}
Write
\[
\Psi_\Delta:=\Phi_{H_0+\Delta H}(\Lambda,0;z),
\qquad
\Psi_0:=\Phi_{H_0}(\Lambda,0;z).
\]
By \Cref{prop:standard}(v),
\[
\Psi_0\in \mathcal V,\qquad \mathfrak m(\Psi_0)=i,\qquad \mathfrak v(\Psi_0)=0.
\]
Since $\mathcal V$ is open, choose $\varepsilon_z>0$ so that the closed ball
\[
K_z:=\{\Psi\in M_2(\C):\|\Psi-\Psi_0\|_{\op}\le \varepsilon_z\}
\]
is contained in $\mathcal V$.

Set
\[
\|\Delta H\|_{L^1}:=\|\Delta H\|_{L^1(0,\Lambda;\op)}.
\]
By \Cref{lem:Phi-diff},
\[
\|\Psi_\Delta-\Psi_0\|_{\op}
\le
|z|e^{2|z|\Lambda}\|\Delta H\|_{L^1}.
\]
Hence, after shrinking $r_v(z,\Lambda)$ if necessary, we may assume that
\[
|z|e^{2|z|\Lambda}r_v(z,\Lambda)\le \varepsilon_z,
\]
so $\Psi_\Delta\in K_z$ whenever $\|\Delta H\|_{L^1}\le r_v(z,\Lambda)$.

Now decompose
\[
\Psi_\Delta-\Psi_0
= D\Phi_{H_0}(\Lambda,0;z)[\Delta H]+E_\Phi,
\]
where \Cref{lem:duhamel} gives
\[
\|E_\Phi\|_{\op}
\le
|z|^2e^{4|z|\Lambda}\|\Delta H\|_{L^1}^2.
\]
Applying \Cref{lem:Phi-to-v} on $K_z$ at the base point $\Psi_0$, we obtain
\begin{align*}
&|\mathfrak v(\Psi_\Delta)-\mathfrak v(\Psi_0) -D\mathfrak v(\Psi_0)[\Psi_\Delta-\Psi_0]| \\
&\qquad\le C_{K_z}\|\Psi_\Delta-\Psi_0\|_{\op}^2
\le
C_{K_z}|z|^2e^{4|z|\Lambda}\|\Delta H\|_{L^1}^2.
\end{align*}
Since $\mathfrak v(\Psi_0)=0$, this yields
\begin{align*}
&|\mathfrak v(\Psi_\Delta)-D\mathfrak v(\Psi_0)[D\Phi_{H_0}(\Lambda,0;z)[\Delta H]]| \\
&\qquad\le |D\mathfrak v(\Psi_0)|\,\|E_\Phi\|_{\op} + C_{K_z}|z|^2e^{4|z|\Lambda}\|\Delta H\|_{L^1}^2.
\end{align*}
The first term on the right is also $O(\|\Delta H\|_{L^1}^2)$ by the bound on $E_\Phi$.

Finally,
\[
\mathfrak v(\Psi_\Delta)=v_{H_0+\Delta H,\Lambda}(z)
\]
by the transfer formula \eqref{eq:weyl-transfer}, while the chain rule gives
\[
D\mathfrak v(\Psi_0)[D\Phi_{H_0}(\Lambda,0;z)[\Delta H]] = Dv_{H_0,\Lambda}(z)[\Delta H].
\]
Combining the bounds and enlarging the constant if necessary proves
\eqref{eq:v-quad}.
\end{proof}

\begin{remark}[Role of the smallness assumption]
The smallness assumption in \Cref{thm:quadratic} is used only to keep the transfer matrix in a fixed neighborhood of the free point on which the passage from transfer matrix to Schur data has a uniform $C^2$ bound.
\end{remark}

\begin{corollary}[Finite-sample free-point expansion]\label{cor:sample-expansion}
Fix sample points $z_k=x_k+i\eta\in\Cplus$. If
\[
\|\Delta H\|_{L^\infty(0,\Lambda;\op)}\le \frac14,
\qquad
\|\Delta H\|_{L^1(0,\Lambda;\op)}\le \min_{1\le k\le M} r_v(z_k,\Lambda),
\]
then
\[
\mathcal S(H_0+\Delta H)
=
\bigl(Dv_{H_0,\Lambda}(z_k)[\Delta H]\bigr)_{k=1}^M
+ R(\Delta H),
\]
with
\[
\|R(\Delta H)\|_2
\le
\sqrt M\,\max_{1\le k\le M} C(z_k,\Lambda)\,
\|\Delta H\|_{L^1(0,\Lambda;\op)}^2.
\]
\end{corollary}

\begin{proof}
Apply \Cref{thm:quadratic} to each sample point $z_k$ and combine the componentwise quadratic bounds with the Euclidean norm estimate.
\end{proof}

\section{Finite-dimensional consequences near the free point}\label{sec:finite}

\subsection{Linear families and the free Jacobian}

Fix $d\ge 1$ and choose $q_1,\dots,q_d\in L^\infty(0,\Lambda;\C)$. For $\theta=(\theta_1,\dots,\theta_d)\in\R^d$, define
\[
q_\theta(s):=\sum_{j=1}^d \theta_j q_j(s),
\qquad
H(\theta)(s):=\frac12 I+\Xmat{q_\theta(s)}.
\]
Set
\[
b(s):=\Bigl(\sum_{j=1}^d |q_j(s)|^2\Bigr)^{1/2},
\qquad
B_\infty:=\|b\|_{L^\infty(0,\Lambda)},
\qquad
B_1:=\|b\|_{L^1(0,\Lambda)}.
\]
If $B_\infty=0$, the family is constant and all subsequent statements are trivial. For a fixed $\varepsilon\in(0,\frac12)$, the family is admissible on
\[
B(0,r_\varepsilon),\qquad r_\varepsilon:=\frac{\frac12-\varepsilon}{B_\infty},
\]
provided $B_\infty>0$.

\begin{proposition}[Free Jacobian formula]\label{prop:general-J}
For $z_k=x_k+i\eta\in\Cplus$,
\[
(D\mathcal S(0)\theta)_k = -iz_k \int_0^\Lambda q_\theta(s)e^{iz_k s}\,ds = -iz_k \sum_{j=1}^d \theta_j \widehat q_j(z_k),
\]
where
\[
\widehat q_j(z):=\int_0^\Lambda q_j(s)e^{izs}\,ds.
\]
\end{proposition}

\begin{proof}
Substitute $q_\theta=\sum_{j=1}^d \theta_j q_j$ into \eqref{eq:dv-free} and exchange the finite sum with the integral.
\end{proof}

\subsection{A smooth dependence lemma}

\begin{lemma}[C2-dependence of the transfer matrix on affine parameters]
\label{lem:C2-transfer}
Let
\[
H(\theta)=H_*+\sum_{j=1}^d \theta_j K_j,
\qquad
\theta\in B(0,r)\subset\R^d,
\]
where $H_*$ is trace-normed on $[0,\Lambda]$, each $K_j\in L^\infty(0,\Lambda;\R^{2\times2})$, and $H(\theta)$ remains nonnegative and trace-normed for $\theta\in B(0,r)$. Fix $z\in\C$. Then
\[
\theta\longmapsto \Phi_{H(\theta)}(\Lambda,0;z)
\]
is a $C^2$ map from $B(0,r)$ into $M_2(\C)$. Its derivatives are given by
\begin{align}
D\Phi_\theta(\Lambda,0;z)[\xi] &= -z\int_0^\Lambda
\Phi_\theta(\Lambda,s;z)\,\Jmat\,K_\xi(s)\,\Phi_\theta(s,0;z)\,ds,
\label{eq:C2-transfer-D1}
\\[0.5em]
D^2\Phi_\theta(\Lambda,0;z)[\xi_1,\xi_2] &= z^2\int_0^\Lambda\int_0^s
\Phi_\theta(\Lambda,s;z)\,\Jmat\,K_{\xi_1}(s)\,\Phi_\theta(s,u;z)
\notag\\
&\qquad\cdot\,\Jmat\,K_{\xi_2}(u)\,\Phi_\theta(u,0;z)\,du\,ds
\notag\\
&\quad+ z^2\int_0^\Lambda\int_0^s
\Phi_\theta(\Lambda,s;z)\,\Jmat\,K_{\xi_2}(s)\,\Phi_\theta(s,u;z)
\notag\\
&\qquad\cdot\,\Jmat\,K_{\xi_1}(u)\,\Phi_\theta(u,0;z)\,du\,ds,
\label{eq:C2-transfer-D2}
\end{align}
where
\[
K_\xi:=\sum_{j=1}^d \xi_j K_j.
\]
\end{lemma}

\begin{proof}
For $\theta\in B(0,r)$, the transfer matrix satisfies the Volterra equation
\[
\Phi_\theta(s,0;z)
= I-z\int_0^s \Jmat H(\theta)(u)\Phi_\theta(u,0;z)\,du.
\]
Since each $H(\theta)$ is nonnegative and trace-normed, the standard propagator bound gives
\[
\|\Phi_\theta(s,t;z)\|_{\op}\le e^{|z|(s-t)}
\qquad
(\theta\in B(0,r),\ 0\le t\le s\le\Lambda).
\]
This bound is uniform on the whole ball.

Let $\theta,\xi\in\R^d$ with $\theta,\theta+\varepsilon\xi\in B(0,r)$. Subtracting the Volterra equations for $\Phi_{\theta+\varepsilon\xi}$ and $\Phi_\theta$ gives
\[
\Phi_{\theta+\varepsilon\xi}(s,0;z)-\Phi_\theta(s,0;z)
= -z\int_0^s
\Phi_{\theta+\varepsilon\xi}(s,u;z)\,\Jmat\,\varepsilon K_\xi(u)\,
\Phi_\theta(u,0;z)\,du.
\]
Dividing by $\varepsilon$, the integrand is dominated by
\[
|z|e^{2|z|\Lambda}\|K_\xi(u)\|_{\op},
\]
which is integrable. Dominated convergence therefore yields the first derivative formula \eqref{eq:C2-transfer-D1}.

To obtain the second derivative, differentiate \eqref{eq:C2-transfer-D1} in the direction $\xi_2$. The derivative may again be passed under the integral sign because each differentiated term is dominated by an integrable majorant of the form
\[
C(z,\Lambda)\,\|K_{\xi_1}(s)\|_{\op}\,\|K_{\xi_2}(u)\|_{\op},
\]
coming from the uniform propagator bound and the first derivative formula already proved. This yields \eqref{eq:C2-transfer-D2}. Continuity of the first and second derivatives follows from the same dominated-convergence argument. Since the parameter dependence is affine, no higher-order parameter derivatives of the Hamiltonian appear.
\end{proof}

\subsection{C2 regularity of the sampling map}

\begin{proposition}[C2 regularity and Jacobian-Lipschitz bound]
\label{prop:C2-family}
Fix sample points $z_k=x_k+i\eta\in\Cplus$, $k=1,\dots,M$, and $\varepsilon\in(0,\frac12)$. Consider the linear family
\[
q_\theta(s):=\sum_{j=1}^d \theta_j q_j(s),
\qquad
H(\theta)(s):=\frac12 I+\Xmat{q_\theta(s)},
\qquad
\theta\in\R^d,
\]
with $q_j\in L^\infty(0,\Lambda;\C)$, and let
\[
b(s):=\Bigl(\sum_{j=1}^d |q_j(s)|^2\Bigr)^{1/2},
\qquad
B_\infty:=\|b\|_{L^\infty(0,\Lambda)},
\qquad
B_1:=\|b\|_{L^1(0,\Lambda)}.
\]
If $B_\infty=0$, the family is constant and there is nothing to prove. Assume $B_\infty>0$, and define
\[
r_\varepsilon:=\frac{\frac12-\varepsilon}{B_\infty}.
\]
Set
\[
Z:=\max_{1\le k\le M}|z_k|,
\qquad
C_1:=ZB_1e^{2Z\Lambda},
\qquad
C_2:=2Z^2B_1^2e^{2Z\Lambda},
\]
and
\[
r_d:=\min\left\{r_\varepsilon,\frac{e^{\eta\Lambda/2}}{4C_1}\right\}.
\]
Then
\[
\mathcal S_\R:B(0,r_d)\subset\R^d\to\R^{2M}
\]
is $C^2$ on $B(0,r_d)$. In particular, there exists
\[
B=B(\Lambda,\eta,z_1,\dots,z_M,q_1,\dots,q_d,\varepsilon)
\]
such that
\[
\|D\mathcal S_\R(\theta)-D\mathcal S_\R(\widetilde\theta)\|_{2\to2}
\le
B\|\theta-\widetilde\theta\|_2,
\qquad
\theta,\widetilde\theta\in B(0,r_d).
\]
\end{proposition}

\begin{proof}
Apply \Cref{lem:C2-transfer} with
\[
H(\theta)=\frac12 I+\Xmat{q_\theta}.
\]
For $\xi\in\R^d$,
\[
q_\xi(s):=\sum_{j=1}^d \xi_j q_j(s),
\qquad
\Delta H_\xi(s):=\Xmat{q_\xi(s)}.
\]
By Cauchy--Schwarz and \Cref{lem:Xnorm},
\[
\|\Delta H_\xi(s)\|_{\op}=|q_\xi(s)|\le b(s)\|\xi\|_2
\]
for a.e.\ $s$, and therefore
\begin{equation}\label{eq:DeltaHxi-bounds}
\|\Delta H_\xi\|_{L^1(0,\Lambda;\op)}\le B_1\|\xi\|_2,
\qquad
\|\Delta H_\xi\|_{L^\infty(0,\Lambda;\op)}\le B_\infty\|\xi\|_2.
\end{equation}
Hence $H(\theta)$ is nonnegative and trace-normed on $B(0,r_\varepsilon)$.

Fix one sample point $z\in\{z_1,\dots,z_M\}$ and write
\[
\Phi_\theta(s,t;z):=\Phi_{H(\theta)}(s,t;z).
\]
From \eqref{eq:C2-transfer-D1} and \eqref{eq:DeltaHxi-bounds},
\begin{align*}
\|D\Phi_\theta(\Lambda,0;z)[\xi]\|_{\op}
&\le |z|\int_0^\Lambda
\|\Phi_\theta(\Lambda,s;z)\|_{\op}\,
\|\Delta H_\xi(s)\|_{\op}\,
\|\Phi_\theta(s,0;z)\|_{\op}\,ds \\
&\le |z|e^{|z|\Lambda}\|\Delta H_\xi\|_{L^1(0,\Lambda;\op)}
\le
|z|B_1e^{|z|\Lambda}\|\xi\|_2
\le C_1\|\xi\|_2.
\end{align*}
Thus
\begin{equation}\label{eq:DPhi-bound}
\|D\Phi_\theta(\Lambda,0;z)\|_{2\to\op}\le C_1.
\end{equation}

Similarly, from \eqref{eq:C2-transfer-D2},
\begin{align*}
\|D^2\Phi_\theta(\Lambda,0;z)[\xi_1,\xi_2]\|_{\op}
&\le 2|z|^2
\int_0^\Lambda\int_0^s
e^{|z|(\Lambda-s)}\|\Delta H_{\xi_1}(s)\|_{\op}\, e^{|z|(s-u)}
\|\Delta H_{\xi_2}(u)\|_{\op}\,
e^{|z|u}\,du\,ds \\
&\le 2|z|^2 e^{|z|\Lambda}
\|\Delta H_{\xi_1}\|_{L^1(0,\Lambda;\op)}
\|\Delta H_{\xi_2}\|_{L^1(0,\Lambda;\op)} \\
&\le 2|z|^2B_1^2e^{|z|\Lambda}\|\xi_1\|_2\|\xi_2\|_2
\le
C_2\|\xi_1\|_2\|\xi_2\|_2.
\end{align*}
Hence
\begin{equation}\label{eq:D2Phi-bound}
\|D^2\Phi_\theta(\Lambda,0;z)\|_{2\to2\to\op}\le C_2.
\end{equation}

By the mean value theorem and \eqref{eq:DPhi-bound},
\[
\|\Phi_\theta(\Lambda,0;z)-\Phi_0(\Lambda,0;z)\|_{\op}
\le
C_1\|\theta\|_2.
\]
Now fix $k\in\{1,\dots,M\}$ and write
\[
\Phi_\theta(\Lambda,0;z_k)=
\begin{pmatrix}
A_k(\theta)&B_k(\theta)\\
C_k(\theta)&D_k(\theta)
\end{pmatrix},
\qquad
d_k(\theta):=D_k(\theta)-iB_k(\theta).
\]
At $\theta=0$,
\[
d_k(0)=e^{-iz_k\Lambda/2},
\qquad
|d_k(0)|=e^{\eta\Lambda/2}
\]
by \Cref{prop:standard}(v). For $\theta\in B(0,r_d)$,
\[
|d_k(\theta)-d_k(0)|
\le
2\|\Phi_\theta(\Lambda,0;z_k)-\Phi_0(\Lambda,0;z_k)\|_{\op}
\le
2C_1\|\theta\|_2
\le
\frac12 e^{\eta\Lambda/2}.
\]
Therefore
\[
|d_k(\theta)|\ge \frac12 e^{\eta\Lambda/2}
\qquad
(\theta\in B(0,r_d)).
\]

Set
\[
n_k(\theta):=iA_k(\theta)-C_k(\theta),
\qquad
m_k(\theta):=\frac{n_k(\theta)}{d_k(\theta)},
\qquad
v_k(\theta):=\frac{m_k(\theta)-i}{m_k(\theta)+i}.
\]
The entries of $\Phi_\theta$ are $C^2$ by \Cref{lem:C2-transfer}, hence so are $n_k$ and $d_k$. Since $d_k$ stays uniformly away from zero on $B(0,r_d)$, the quotient rule shows that $m_k$ is $C^2$ on $B(0,r_d)$, with bounds depending only on the quantities listed in the statement. Because each $m_k(\theta)\in\Cplus$, the Cayley map
\[
c(w)=\frac{w-i}{w+i}
\]
and its derivatives remain uniformly bounded on $\{m_k(\theta):\theta\in B(0,r_d)\}$. Therefore $v_k=c\circ m_k$ is $C^2$ on $B(0,r_d)$, with a uniform bound on $D^2v_k$.

Finally,
\[
\mathcal S_\R(\theta)
=
\bigl(\Re v_1(\theta),\Im v_1(\theta),\dots,\Re v_M(\theta),\Im v_M(\theta)\bigr).
\]
Hence $\mathcal S_\R$ is $C^2$ on $B(0,r_d)$ and
\[
\sup_{\theta\in B(0,r_d)}
\|D^2\mathcal S_\R(\theta)\|_{2\to2\to2}<\infty.
\]
The Jacobian-Lipschitz estimate follows from the mean value theorem in Banach spaces.
\end{proof}

\subsection{Local inverse results}

We record the following standard consequence of an injective base Jacobian and a locally Lipschitz derivative.

\begin{theorem}[Abstract local bi-Lipschitz lemma]\label{thm:local-bl}
Let $\mathcal S_\R$ be the realification of the sampling map for a finite-dimensional free-tail family on a ball $B(0,r_d)\subset\R^d$. Assume:
\begin{enumerate}[label=\textup{(\alph*)}]
\item $\mathcal S_\R$ is $C^1$ on $B(0,r_d)$;
\item there exists $B\ge 0$ such that
\[
\|D\mathcal S_\R(\theta)-D\mathcal S_\R(\widetilde\theta)\|_{2\to2}
\le
B\|\theta-\widetilde\theta\|_2
\qquad
(\theta,\widetilde\theta\in B(0,r_d));
\]
\item the realified Jacobian $J_\R:=D\mathcal S_\R(0)$ satisfies
\[
\alpha:=\smin(J_\R)>0.
\]
\end{enumerate}
Set
\[
r_*:=\min\left\{r_d,\frac{\alpha}{2B}\right\},
\]
with the convention $\alpha/(2B)=+\infty$ when $B=0$. Then for all $\theta,\widetilde\theta\in B(0,r_*)$,
\[
\frac{\alpha}{2}\|\theta-\widetilde\theta\|_2
\le
\|\mathcal S(\theta)-\mathcal S(\widetilde\theta)\|_2
\le
\bigl(\|J_\R\|_{2\to2}+Br_*\bigr)\|\theta-\widetilde\theta\|_2.
\]
In particular, $\mathcal S$ is injective on $B(0,r_*)$ with respect to the Euclidean parameter metric.
\end{theorem}

\begin{proof}
For $\theta,\widetilde\theta\in B(0,r_*)$, the segment joining them lies in $B(0,r_*)$. Hence
\[
\mathcal S_\R(\theta)-\mathcal S_\R(\widetilde\theta)
=
\int_0^1
D\mathcal S_\R(\widetilde\theta+t(\theta-\widetilde\theta)) (\theta-\widetilde\theta)\,dt.
\]
Subtract and add $J_\R(\theta-\widetilde\theta)$ and use the Jacobian-Lipschitz bound:
\begin{align*}
\|\mathcal S(\theta)-\mathcal S(\widetilde\theta)\|_2
&\ge
\|J_\R(\theta-\widetilde\theta)\|_2
-
\int_0^1
\|D\mathcal S_\R(\widetilde\theta+t(\theta-\widetilde\theta))-J_\R\|_{2\to2}
\,dt\,
\|\theta-\widetilde\theta\|_2 \\
&\ge (\alpha-Br_*)\|\theta-\widetilde\theta\|_2
\ge
\frac{\alpha}{2}\|\theta-\widetilde\theta\|_2.
\end{align*}
The upper bound is proved in the same way.
\end{proof}

\begin{remark}
The only model-dependent quantity in \Cref{thm:local-bl} is the least singular value of the free Jacobian. For the linear families considered here, the nonlinear bounds come uniformly from the Volterra equation through \Cref{prop:C2-family}.
\end{remark}

In the square case the preceding lemma yields the usual frozen-Jacobian contraction scheme.

\begin{theorem}[Abstract local inverse chart in the square case]\label{thm:IFT}
Assume the square case $d=2M$ and the hypotheses of \Cref{thm:local-bl}. Set
\[
M_0:=\alpha^{-1},
\qquad
r:=\min\left\{r_d,\frac{1}{2BM_0}\right\},
\qquad
\delta:=\frac{r}{2M_0},
\]
with the convention $1/(2BM_0)=+\infty$ when $B=0$. Then:
\begin{enumerate}[label=\textup{(\roman*)}]
\item for every $y\in\C^M$ with $\|y-\mathcal S(0)\|_2\le \delta$, there exists a
unique $\theta\in B(0,r)$ such that $\mathcal S(\theta)=y$;

\item the local inverse is Lipschitz:
\[
\|\mathcal S^{-1}(y)-\mathcal S^{-1}(\widetilde y)\|_2
\le
2M_0\|y-\widetilde y\|_2;
\]

\item the map
\[
\mathcal G_{y,\R}(\theta):=
\theta-J_\R^{-1}\bigl(\mathcal S_\R(\theta)-y_\R\bigr)
\]
is a contraction on $B(0,r)$ and converges to the unique local solution.
\end{enumerate}
\end{theorem}

\begin{proof}
Set
\[
F_\R(\theta):=\mathcal S_\R(\theta)-y_\R,
\qquad
\mathcal G_{y,\R}(\theta):=\theta-J_\R^{-1}F_\R(\theta).
\]
Since $M_0=\|J_\R^{-1}\|_{2\to2}=\alpha^{-1}$,
\[
D\mathcal G_{y,\R}(\theta) = I-J_\R^{-1}D\mathcal S_\R(\theta).
\]
Hence, for $\theta\in B(0,r)$,
\[
\|D\mathcal G_{y,\R}(\theta)\|_{2\to2}
\le
M_0\|D\mathcal S_\R(\theta)-J_\R\|_{2\to2}
\le
M_0Br
\le \frac12.
\]
Therefore $\mathcal G_{y,\R}$ is $\frac12$-Lipschitz on $B(0,r)$.

Moreover,
\[
\mathcal G_{y,\R}(0)=J_\R^{-1}(y_\R-\mathcal S_\R(0)),
\]
so
\[
\|\mathcal G_{y,\R}(0)\|_2
\le
M_0\|y_\R-\mathcal S_\R(0)\|_2
\le
M_0\delta =\frac r2.
\]
Thus, for every $\theta\in B(0,r)$,
\[
\|\mathcal G_{y,\R}(\theta)\|_2
\le
\|\mathcal G_{y,\R}(\theta)-\mathcal G_{y,\R}(0)\|_2+\|\mathcal G_{y,\R}(0)\|_2
\le
\frac12\|\theta\|_2+\frac r2
\le r.
\]
Hence $\mathcal G_{y,\R}$ maps $B(0,r)$ into itself and is a contraction there. Banach's fixed-point theorem gives the unique solution in $B(0,r)$ and the convergence of the iteration.

For the Lipschitz bound on the inverse, let
\[
\theta:=\mathcal S^{-1}(y),\qquad \widetilde\theta:=\mathcal S^{-1}(\widetilde y).
\]
Since $\theta,\widetilde\theta\in B(0,r)\subset B(0,r_*)$, \Cref{thm:local-bl} yields
\[
\frac{\alpha}{2}\|\theta-\widetilde\theta\|_2
\le
\|\mathcal S(\theta)-\mathcal S(\widetilde\theta)\|_2
=
\|y-\widetilde y\|_2.
\]
This gives
\[
\|\mathcal S^{-1}(y)-\mathcal S^{-1}(\widetilde y)\|_2
\le
2\alpha^{-1}\|y-\widetilde y\|_2 = 2M_0\|y-\widetilde y\|_2.
\]
\end{proof}

\section{The block model: exact factorization and depth bounds}\label{sec:block}

\subsection{Block Hamiltonians}

Fix $N\ge 1$ and set $\ell:=\Lambda/N$. Partition
\[
[0,\Lambda]=\bigcup_{j=0}^{N-1} I_j,
\qquad
I_j=[j\ell,(j+1)\ell).
\]
A block Hamiltonian is constant on each $I_j$:
\[
H(s)=H_j\quad \text{for } s\in I_j,
\]
with
\[
H_j=\frac12 I+\Xmat{q_j},
\qquad
q_j=h_j+ik_j\in\C.
\]
The parameter is $q=(q_0,\dots,q_{N-1})\in\C^N$.

\subsection{Exact factorization}

For a fixed design $x=(x_1,\dots,x_M)\in\R^M$, let
\[
T_x:\C^N\to\C^M
\]
denote the complex-linear form of the free Jacobian in block coordinates.

\begin{theorem}[Exact block linearization]\label{thm:block-factor}
Fix $z=x+i\eta\in\Cplus$. For a block perturbation $\Delta q\in\C^N$,
\begin{equation}\label{eq:block-derivative}
Dv_{H_0,\Lambda}(z)[\Delta q] = -2i\sin\!\Bigl(\frac{z\ell}{2}\Bigr)
\sum_{j=0}^{N-1}\Delta q_j\,e^{iz(j+\tfrac12)\ell}.
\end{equation}
Hence, for the finite sample set $z_k=x_k+i\eta$,
\begin{equation}\label{eq:T-factor}
(T_x\Delta q)_k =
\gamma_k(x_k)\sum_{j=0}^{N-1}\Delta q_j\,e^{ix_kj\ell}w_j,
\end{equation}
where
\[
\gamma_k(x_k):=
-2i\sin\!\Bigl(\frac{(x_k+i\eta)\ell}{2}\Bigr)e^{ix_k\ell/2},
\qquad
w_j:=e^{-\eta(j+\tfrac12)\ell}.
\]
Equivalently,
\[
T_x=D_\gamma(x)\,F_x\,D_w,
\]
with
\[
(F_x)_{kj}=e^{ix_kj\ell},
\qquad
D_\gamma(x)=\diag(\gamma_1(x_1),\dots,\gamma_M(x_M)),
\qquad
D_w=\diag(w_0,\dots,w_{N-1}).
\]
\end{theorem}

\begin{proof}
The block model reduces the free derivative to exact integration of exponentials over the block intervals. Applying \eqref{eq:dv-free} with $q(s)\equiv \Delta q_j$ on $I_j$ gives
\begin{align*}
Dv_{H_0,\Lambda}(z)[\Delta q] &= -iz\sum_{j=0}^{N-1}\Delta q_j \int_{j\ell}^{(j+1)\ell} e^{izs}\,ds \\
&= -\sum_{j=0}^{N-1}\Delta q_j\bigl(e^{iz(j+1)\ell}-e^{izj\ell}\bigr) \\
&= -2i\sin\!\Bigl(\frac{z\ell}{2}\Bigr)
\sum_{j=0}^{N-1}\Delta q_j\,e^{iz(j+\tfrac12)\ell}.
\end{align*}
Now split
\[
e^{iz(j+\tfrac12)\ell} = e^{ix\ell/2}\,e^{ixj\ell}\,e^{-\eta(j+\tfrac12)\ell},
\]
which yields \eqref{eq:T-factor}.
\end{proof}

\begin{remark}
This factorization is exact at the free point and isolates row geometry, Fourier geometry, and exponential depth attenuation. It shows in particular that depth ill-conditioning is already present at first order and is not an artifact of the nonlinear inverse problem.
\end{remark}

\subsection{The Fourier factor}

Define the Fourier sampling matrix
\[
(F_x)_{kj}:=e^{ix_kj\ell},
\qquad 1\le k\le M,\ 0\le j\le N-1.
\]

\begin{lemma}[Singular values under factors]\label{lem:smin-prod}
For compatible matrices $A$, $B$, $C$,
\[
\smin(ABC)\ge \smin(A)\smin(B)\smin(C).
\]
\end{lemma}

\begin{proof}
For every unit vector $u$,
\[
\|ABC\,u\|_2
\ge
\smin(A)\|BCu\|_2
\ge
\smin(A)\smin(B)\|Cu\|_2
\ge
\smin(A)\smin(B)\smin(C).
\]
\end{proof}

\begin{lemma}[Tight-frame bound for the Fourier factor]\label{lem:tight}
If $M\ge N$, then
\[
\smin(F_x)\le \sqrt M.
\]
Equality holds if and only if
\[
F_x^*F_x=MI_N.
\]
\end{lemma}

\begin{proof}
If $\sigma_1\ge\cdots\ge\sigma_N$ are the singular values of $F_x$, then
\[
\sum_{j=1}^N \sigma_j^2=\|F_x\|_F^2=MN.
\]
Hence $\sigma_N^2\le M$. Equality holds exactly when all singular values are equal to $\sqrt M$, equivalently when $F_x^*F_x=MI_N$.
\end{proof}

\begin{lemma}[Equispaced designs are Fourier-tight]\label{lem:equi-tight}
Assume $M\ge N$. For
\[
x_k(\alpha):=\frac{2\pi(k-1)}{M\ell}+\frac{2\alpha}{\ell},
\qquad k=1,\dots,M,
\]
one has
\[
F_{x(\alpha)}^*F_{x(\alpha)}=MI_N.
\]
\end{lemma}

\begin{proof}
For $0\le j,j'\le N-1$,
\[
(F^*F)_{jj'} =
\sum_{k=1}^M e^{-ix_k(\alpha)j\ell}e^{ix_k(\alpha)j'\ell}
=
\sum_{k=1}^M e^{2\pi i(k-1)(j'-j)/M}.
\]
This equals $M$ when $j=j'$ and $0$ otherwise, because in the off-diagonal case
\[
1\le |j'-j|\le N-1\le M-1.
\]
\end{proof}

\begin{lemma}[Row-factor envelope]\label{lem:row}
For $a\in\R$ and $b>0$,
\[
|\sin(a+ib)|^2=\sin^2 a+\sinh^2 b.
\]
Hence
\[
\sinh b\le |\sin(a+ib)|\le \cosh b.
\]
\end{lemma}

\begin{proof}
Use
\[
\sin(a+ib)=\sin a\cosh b+i\cos a\sinh b.
\]
\end{proof}

\subsection{Bounds in the block model}

\begin{theorem}[Universal upper bound]\label{thm:upper}
For every design $x\in\R^M$,
\[
\smin(T_x)
\le
2\sqrt M\,\cosh\!\Bigl(\frac{\eta\ell}{2}\Bigr)e^{-\eta(\Lambda-\ell/2)}.
\]
\end{theorem}

\begin{proof}
Let $e_{N-1}\in\C^N$ be the deepest-block coordinate vector. Then
\[
\smin(T_x)\le \|T_xe_{N-1}\|_2.
\]
By \eqref{eq:T-factor},
\[
(T_xe_{N-1})_k=\gamma_k(x_k)e^{ix_k(N-1)\ell}w_{N-1},
\]
hence
\[
\|T_xe_{N-1}\|_2^2
= w_{N-1}^2\sum_{k=1}^M |\gamma_k(x_k)|^2.
\]
Now
\[
w_{N-1}=e^{-\eta(\Lambda-\ell/2)},
\]
and by \Cref{lem:row},
\[
|\gamma_k(x_k)| = 2\left|\sin\!\Bigl(\frac{(x_k+i\eta)\ell}{2}\Bigr)\right|
\le
2\cosh\!\Bigl(\frac{\eta\ell}{2}\Bigr).
\]
Taking square roots proves the claim.
\end{proof}

\begin{theorem}[Lower bound for the half-shifted design]\label{thm:lower}
Assume $M\ge N$. For
\[
x_k:=\frac{2\pi(k-\tfrac12)}{M\ell},
\qquad
k=1,\dots,M,
\]
one has
\[
\smin(T_x)
\ge
2\sqrt M\,
\sqrt{\sinh^2\!\Bigl(\frac{\eta\ell}{2}\Bigr)+\sin^2\!\Bigl(\frac{\pi}{2M}\Bigr)}
\,e^{-\eta(\Lambda-\ell/2)}.
\]
\end{theorem}

\begin{proof}
By \Cref{lem:equi-tight},
\[
\smin(F_x)=\sqrt M.
\]
By \Cref{lem:smin-prod},
\[
\smin(T_x)\ge \smin(D_\gamma)\smin(F_x)\smin(D_w).
\]
Now
\[
\smin(D_w)=\min_j w_j=e^{-\eta(\Lambda-\ell/2)}.
\]
Also
\[
\frac{x_k\ell}{2}=\frac{\pi(k-\tfrac12)}{M},
\]
so by \Cref{lem:row},
\[
|\gamma_k(x_k)|^2 = 4\left(
\sinh^2\!\Bigl(\frac{\eta\ell}{2}\Bigr)
+
\sin^2\!\Bigl(\frac{\pi(k-\tfrac12)}{M}\Bigr)
\right).
\]
The minimum occurs at $k=1$ or $k=M$ and equals
\[
4\left(
\sinh^2\!\Bigl(\frac{\eta\ell}{2}\Bigr)
+
\sin^2\!\Bigl(\frac{\pi}{2M}\Bigr)
\right).
\]
Combining the three factors gives the result.
\end{proof}

\begin{theorem}[Best shift for the explicit product lower bound]\label{thm:shift}
Fix $M\ge N$ and consider the equispaced Fourier-tight family
\[
x_k(\alpha)=\frac{2\pi(k-1)}{M\ell}+\frac{2\alpha}{\ell},
\qquad k=1,\dots,M.
\]
Then
\[
\max_{\alpha\in\R}\min_{1\le k\le M}|\gamma_k(x_k(\alpha))|
= 2\sqrt{
\sinh^2\!\Bigl(\frac{\eta\ell}{2}\Bigr)
+
\sin^2\!\Bigl(\frac{\pi}{2M}\Bigr)
}.
\]
The maximizers are exactly
\[
\alpha\equiv \frac{\pi}{2M}\pmod{\frac{\pi}{M}}.
\]
Equivalently, this shift maximizes the lower bound obtained from
\[
\smin(T_x)\ge \smin(D_\gamma)\smin(F_x)\smin(D_w)
\]
within this family.
\end{theorem}

\begin{proof}
By \Cref{lem:equi-tight}, every design in this family satisfies
\[
\smin(F_{x(\alpha)})=\sqrt M.
\]
Hence optimizing the product lower bound from \Cref{lem:smin-prod} reduces to optimizing the worst row factor
\[
\min_{1\le k\le M} |\gamma_k(x_k(\alpha))|.
\]
Now
\[
\frac{x_k(\alpha)\ell}{2}
=
\alpha+\frac{\pi(k-1)}{M},
\]
so by \Cref{lem:row},
\[
|\gamma_k(x_k(\alpha))|^2 = 4\left(
\sinh^2\!\Bigl(\frac{\eta\ell}{2}\Bigr)
+
\sin^2\!\Bigl(\alpha+\frac{\pi(k-1)}{M}\Bigr)
\right).
\]
Thus the problem is equivalent to maximizing
\[
m(\alpha):=
\min_{1\le k\le M}
\sin^2\!\Bigl(\alpha+\frac{\pi(k-1)}{M}\Bigr).
\]

Consider the equally spaced grid
\[
G_\alpha:=
\left\{
\alpha+\frac{\pi(k-1)}{M}\ \bmod \pi:\ k=1,\dots,M
\right\}
\subset \R/\pi\Z.
\]
If
\[
d(\alpha):=\dist_{\R/\pi\Z}(0,G_\alpha),
\]
then
\[
m(\alpha)=\sin^2(d(\alpha)),
\]
because the zeros of $\sin$ in $\R/\pi\Z$ are exactly the class of $0$. Since the mesh is $\pi/M$, the largest possible value of $d(\alpha)$ is $\pi/(2M)$, attained exactly when $0$ is the midpoint between two consecutive grid points. This happens precisely when
\[
\alpha\equiv \frac{\pi}{2M}\pmod{\frac{\pi}{M}}.
\]
Substituting back proves the claim.
\end{proof}

\begin{remark}
\Cref{thm:shift} identifies the optimizer of the worst row factor, and therefore of
the explicit product lower bound coming from the factorization $T_x=D_\gamma(x)F_xD_w$. It does not claim that the same shift maximizes the exact quantity $\smin(T_x)$ over the whole equispaced Fourier-tight family.
\end{remark}

\begin{corollary}[Two-sided design bound]\label{cor:sandwich}
For
\[
\Sigma_M^*(\eta;\Lambda,N):=\sup_{x\in\R^M}\smin(T_x),
\]
one has
\[
2\sqrt M\,
\sqrt{\sinh^2\!\Bigl(\frac{\eta\ell}{2}\Bigr)+\sin^2\!\Bigl(\frac{\pi}{2M}\Bigr)}
\,e^{-\eta(\Lambda-\ell/2)}
\le
\Sigma_M^*(\eta;\Lambda,N)
\le
2\sqrt M\,\cosh\!\Bigl(\frac{\eta\ell}{2}\Bigr)e^{-\eta(\Lambda-\ell/2)}.
\]
In particular, $\Sigma_M^*(\eta;\Lambda,N)$ is trapped between two explicit multiples of $\sqrt M\,e^{-\eta(\Lambda-\ell/2)}$.
\end{corollary}

\begin{proof}
The lower bound comes from \Cref{thm:lower}, and the upper bound is
\Cref{thm:upper}.
\end{proof}

\subsection{Local inversion and a deterministic worst-case lower bound in the block model}

\begin{lemma}[Realification preserves singular values]\label{lem:realification}
Let $T:\C^N\to\C^M$ be complex-linear, and let $T_\R:\R^{2N}\to\R^{2M}$ be its realification. Then the singular values of $T_\R$ are exactly the singular values of $T$, each repeated twice. In particular,
\[
\smin(T_\R)=\smin(T).
\]
\end{lemma}

\begin{proof}
If $T=U\Sigma V^*$ is a singular value decomposition over $\C$, then the realifications $U_\R$ and $V_\R$ are orthogonal and
\[
T_\R=U_\R \Sigma_\R V_\R^\top,
\]
where $\Sigma_\R$ is block diagonal and each singular value of $T$ appears twice.
\end{proof}

\begin{corollary}[Local inversion in the square block model]\label{cor:block-IFT}
Assume $M=N\ge 2$ and choose the half-shifted equispaced design. Let $\mathcal S_\R$ denote the realified sampling map in the $2N$ real block coordinates, and let $B$ be the Jacobian-Lipschitz constant furnished by
\Cref{prop:C2-family} for the block family. Then the conclusions of
\Cref{thm:IFT} hold with
\[
\alpha
\ge
2\sqrt N\,
\sqrt{\sinh^2\!\Bigl(\frac{\eta\ell}{2}\Bigr)+\sin^2\!\Bigl(\frac{\pi}{2N}\Bigr)}
\,e^{-\eta(\Lambda-\ell/2)}.
\]
In particular, the local inverse has Lipschitz constant at most
\[
\frac{
e^{\eta(\Lambda-\ell/2)} }{
\sqrt N\,
\sqrt{\sinh^2(\eta\ell/2)+\sin^2(\pi/(2N))}
}.
\]
Here the lower bound on $\alpha$ comes entirely from the explicit free Jacobian factorization and the realification lemma.
\end{corollary}

\begin{proof}
The block family is a special case of \Cref{prop:C2-family}. By
\Cref{lem:realification}, the realified Jacobian has the same least singular value as
the complex Jacobian, and \Cref{thm:lower} gives the lower bound on $\alpha$.
\end{proof}

\begin{theorem}[A deterministic worst-case lower bound in the block model]
\label{thm:block-minimax}
Assume the square block model $M=N$ and let $\alpha:=\smin(D\mathcal S_\R(0))>0$. Let $B$ be the Jacobian-Lipschitz constant furnished by
\Cref{prop:C2-family} on a base ball $B(0,r_d)$. In the deterministic bounded-noise
model with noise radius $\delta>0$, set
\[
t_0:=\min\left\{r_d,\frac{\delta}{2\alpha}\right\},
\qquad
t:=
\begin{cases}
\min\{t_0,\sqrt{\delta/B}\}, & B>0,\\
t_0, & B=0.
\end{cases}
\]
Then in the observation model
\[
y=\mathcal S(\theta)+e,\qquad \|e\|_2\le \delta,\qquad \theta\in B(0,r_d),
\]
every reconstruction rule $\widehat\theta$ satisfies
\[
\sup_{\|e\|_2\le \delta}\ \sup_{\theta\in B(0,r_d)}
\|\widehat\theta-\theta\|_2
\ge \frac{t}{2}.
\]
If
\[
\delta\le 2\alpha r_d
\quad\text{and, if } B>0,\quad
\delta\le \frac{\alpha^2}{B},
\]
then
\[
\sup_{\|e\|_2\le \delta}\ \sup_{\theta\in B(0,r_d)}
\|\widehat\theta-\theta\|_2
\ge \frac{\delta}{4\alpha}.
\]
\end{theorem}

\begin{proof}
Choose a unit right-singular vector $h$ of $D\mathcal S_\R(0)$ such that
\[
\|D\mathcal S_\R(0)h\|_2=\alpha.
\]
Let $\theta_0=0$ and $\theta_1=th$. A second-order Taylor estimate gives
\[
\|\mathcal S_\R(\theta_1)-\mathcal S_\R(\theta_0)\|_2
\le
t\alpha+\frac B2 t^2.
\]
By the choice of $t$, the right-hand side is at most $\delta$.

Now consider the single observation
\[
y:=\mathcal S(\theta_0).
\]
Under $\theta_0$, this corresponds to noise $e_0=0$. Under $\theta_1$, it corresponds to
\[
e_1:=\mathcal S(\theta_0)-\mathcal S(\theta_1),
\qquad
\|e_1\|_2\le \delta.
\]
Thus the same observation is compatible with both $\theta_0$ and $\theta_1$ under admissible noise, so no reconstruction rule can guarantee error smaller than $t/2$ uniformly over the admissible parameter-noise pairs.
\end{proof}

\section{Finite-sample invisibility on the full free-tail class}\label{sec:invis}

\begin{theorem}[Deep first-order invisible directions]\label{thm:invis}
Fix $\eta>0$ and sample points
\[
z_k=x_k+i\eta\in\Cplus,\qquad k=1,\dots,M.
\]
Fix also a depth threshold $s_\star\in[0,\Lambda)$. Then there exists a nonzero step function
\[
q\in L^\infty(0,\Lambda;\C),
\qquad
\supp q\subset [s_\star,\Lambda],
\]
such that
\begin{equation}\label{eq:kernel}
Dv_{H_0,\Lambda}(z_k)[q]=0,\qquad k=1,\dots,M.
\end{equation}
Thus every finite fixed-height sample set annihilates nontrivial perturbation directions supported arbitrarily close to the seam. After rescaling we may assume $\|q\|_{L^\infty}\le 1$.

Let
\[
H_\pm^\tau(s):=\frac12 I\pm \tau \Xmat{q(s)},
\qquad 0\le s\le \Lambda,
\]
and extend by the free tail for $s\ge \Lambda$. Then for every
\[
0<\tau\le \tau_0:=
\min\left\{
\frac14,\,
\min_{1\le k\le M}\frac{r_v(z_k,\Lambda)}{\|q\|_{L^1(0,\Lambda)}}
\right\},
\]
the Hamiltonians $H_\pm^\tau$ are admissible and satisfy
\begin{equation}\label{eq:seam-gap}
\|\mathcal S(H_+^\tau)-\mathcal S(H_-^\tau)\|_2
\le
2\sqrt M\,C_*\,\tau^2\,\|q\|_{L^1(0,\Lambda)}^2,
\end{equation}
where
\[
C_*:=\max_{1\le k\le M} C(z_k,\Lambda).
\]
At the same time,
\begin{equation}\label{eq:L1-gap}
\|H_+^\tau-H_-^\tau\|_{L^1(0,\Lambda;\op)}
= 2\tau\,\|q\|_{L^1(0,\Lambda)}.
\end{equation}
\end{theorem}

\begin{proof}
Choose $M+1$ pairwise disjoint subintervals
\[
J_0,\dots,J_M\subset [s_\star,\Lambda].
\]
For each $j$, define
\[
u_j:=\left(\int_{J_j} e^{iz_k s}\,ds\right)_{k=1}^M \in \C^M.
\]
The vectors $u_0,\dots,u_M\in\C^M$ are linearly dependent because their number exceeds the ambient dimension. Hence there exist coefficients $a_0,\dots,a_M$, not all zero, such that
\[
\sum_{j=0}^M a_j u_j=0.
\]
Set
\[
q(s):=\sum_{j=0}^M a_j\mathbf 1_{J_j}(s).
\]
Then $q\neq 0$, $\supp q\subset [s_\star,\Lambda]$, and
\[
\int_0^\Lambda q(s)e^{iz_k s}\,ds=0,\qquad k=1,\dots,M.
\]
By \eqref{eq:dv-free}, this is exactly \eqref{eq:kernel}.

After rescaling, $\|q\|_{L^\infty}\le 1$. By \Cref{lem:Xnorm},
\[
\|\tau\Xmat{q(s)}\|_{\op}\le \tau\le \frac14.
\]
Hence \Cref{lem:free-admissible} shows that $H_\pm^\tau$ are admissible trace-normed Hamiltonians.

Fix $k$. Apply \Cref{thm:quadratic} at $z=z_k$ to the perturbations
\[
\Delta H_\pm=\pm\tau\Xmat{q}.
\]
Because the linear term vanishes by \eqref{eq:kernel},
\[
|v_{H_\pm^\tau,\Lambda}(z_k)|
\le
C(z_k,\Lambda)\tau^2\|q\|_{L^1(0,\Lambda)}^2.
\]
Therefore
\[
|v_{H_+^\tau,\Lambda}(z_k)-v_{H_-^\tau,\Lambda}(z_k)|
\le
2C_* \tau^2\|q\|_{L^1(0,\Lambda)}^2.
\]
Summing over $k$ gives \eqref{eq:seam-gap}. Finally,
\[
H_+^\tau-H_-^\tau=2\tau \Xmat{q},
\]
so \eqref{eq:L1-gap} follows from \Cref{lem:Xnorm}.
\end{proof}

\begin{remark}
This obstruction comes from finite codimension of the sampled Fourier--Laplace data at the free point, not from any special feature of step functions.
\end{remark}

\begin{corollary}[Failure of local inverse-Lipschitz stability]\label{cor:no-lip}
For the sampling map on the full free-tail class,
\[
\frac{
\|\mathcal S(H_+^\tau)-\mathcal S(H_-^\tau)\|_2
}{
\|H_+^\tau-H_-^\tau\|_{L^1(0,\Lambda;\op)}
}
\longrightarrow 0
\qquad
(\tau\downarrow 0).
\]
Hence there is no $L^1(0,\Lambda;\op)$-neighborhood of $H_0$ on which a lower Lipschitz estimate
\[
\|\mathcal S(H_1)-\mathcal S(H_2)\|_2
\ge c\,\|H_1-H_2\|_{L^1(0,\Lambda;\op)}
\]
holds for all free-tail Hamiltonians $H_1,H_2$.
\end{corollary}

\begin{proof}
Divide \eqref{eq:seam-gap} by \eqref{eq:L1-gap}. This obstruction is stated for the full free-tail class endowed with the Hamiltonian $L^1$-distance and does not contradict the finite-dimensional local inversion results proved earlier.
\end{proof}

\begin{corollary}[Deterministic worst-case lower bound on the full class]
\label{cor:minimax-free}
Fix $R>0$ and let
\[
\mathfrak H_R:=
\left\{
H:\ H \text{ is free-tail and }
\|H-H_0\|_{L^1(0,\Lambda;\op)}\le R
\right\}.
\]
In the deterministic bounded-noise observation model
\[
y=\mathcal S(H)+e,\qquad \|e\|_2\le \delta,\qquad H\in\mathfrak H_R,
\]
there exists $c>0$ such that for all sufficiently small $\delta$,
\[
\inf_{\widehat H}\ \sup_{\|e\|_2\le \delta}\ \sup_{H\in\mathfrak H_R}
\|\widehat H-H\|_{L^1(0,\Lambda;\op)}
\ge c\,\sqrt{\delta}.
\]
\end{corollary}

\begin{proof}
Fix a nonzero step function $q$ furnished by \Cref{thm:invis}, and set
\[
c_0:=2\sqrt M\,C_*\,\|q\|_{L^1(0,\Lambda)}^2,
\qquad
c_1:=2\|q\|_{L^1(0,\Lambda)}.
\]
Then \Cref{thm:invis} gives
\[
\|\mathcal S(H_+^\tau)-\mathcal S(H_-^\tau)\|_2\le c_0\tau^2,
\qquad
\|H_+^\tau-H_-^\tau\|_{L^1(0,\Lambda;\op)}=c_1\tau.
\]
Choose
\[
\tau:=\min\left\{\tau_0,\frac{R}{\|q\|_{L^1(0,\Lambda)}},\sqrt{\delta/c_0}\right\}.
\]
Then $H_\pm^\tau\in\mathfrak H_R$. For sufficiently small $\delta$,
\[
\tau=\sqrt{\delta/c_0},
\qquad
\|\mathcal S(H_+^\tau)-\mathcal S(H_-^\tau)\|_2\le \delta.
\]
The same observation is compatible with both $H_+^\tau$ and $H_-^\tau$ under admissible noise of size at most $\delta$, so any estimator must incur at least half their separation for one of the two parameters:
\[
\frac12\|H_+^\tau-H_-^\tau\|_{L^1(0,\Lambda;\op)}
=
\frac{c_1}{2}\tau
=
\frac{c_1}{2\sqrt{c_0}}\sqrt{\delta}.
\]
\end{proof}

\appendix
\section{The normalized Paley--Wiener model}\label{app:pw}

This appendix records only the linearized free-point identification with the normalized Paley--Wiener model; none of the nonlinear results in the main text depends on this reformulation.

Set $a:=\Lambda/2$. For $q\in L^2(0,\Lambda)$ define
\[
(\mathcal U q)(z):=(2\pi)^{-1/2}e^{-iaz}\int_0^\Lambda q(s)e^{izs}\,ds.
\]
Then $\mathcal U$ is unitary from $L^2(0,\Lambda)$ onto the standard Paley--Wiener space $PW_a$; see
\cite{deBranges1968,DymMcKean1976,Remling2018,Romanov2014}. Combining this with
\eqref{eq:dv-free}, we obtain
\[
Dv_{H_0,\Lambda}(z)[q] = -iz\,(2\pi)^{1/2}e^{iaz}(\mathcal U q)(z).
\]
Thus, after division by the explicit multiplier $-iz\,e^{iaz}$, the normalized linearized sampling operator is the usual vertical evaluation map
\[
q\longmapsto (\mathcal U q)(x+i\eta)
\]
for $PW_a$.

Under the Fourier-side model of $PW_a$, evaluation on the horizontal line $\Im z=\eta$ differs from evaluation on the real line by multiplication of the band-limited Fourier data by the factor $e^{-\eta\xi}$ on $\xi\in[-a,a]$. Since this multiplier is bounded above and below on the compact band $[-a,a]$, stable sampling on $\Im z=\eta$ is equivalent to stable sampling on the real line. Accordingly, the normalized free linearization is governed by the classical theory of Paley--Wiener sampling and Fourier frames, in particular the work of Ortega--Cerd\`a and Seip
\cite{OrtegaCerdaSeip2002}.


\begin{thebibliography}{99}

\bibitem{Bessonov2018}
R.~Bessonov,
\emph{Sampling measures, Muckenhoupt Hamiltonians, and triangular factorization},
Int. Math. Res. Not. IMRN 2018, no.~12, 3744--3768.

\bibitem{BessonovRomanov2016}
R.~V.~Bessonov and R.~V.~Romanov,
\emph{An inverse problem for weighted Paley--Wiener spaces},
Inverse Problems \textbf{32} (2016), no.~11, 115007.

\bibitem{Christensen2016}
O.~Christensen,
\emph{An Introduction to Frames and Riesz Bases},
2nd ed., Applied and Numerical Harmonic Analysis, Springer, 2016.

\bibitem{CoddingtonLevinson1955}
E.~A.~Coddington and N.~Levinson,
\emph{Theory of Ordinary Differential Equations},
McGraw--Hill, New York, 1955.

\bibitem{deBranges1968}
L.~de~Branges,
\emph{Hilbert Spaces of Entire Functions},
Prentice--Hall, Englewood Cliffs, NJ, 1968.

\bibitem{DymMcKean1976}
H.~Dym and H.~P.~McKean,
\emph{Gaussian Processes, Function Theory, and the Inverse Spectral Problem},
Academic Press, New York, 1976.

\bibitem{EKT2018}
J.~Eckhardt, A.~Kostenko, and G.~Teschl,
\emph{Spectral asymptotics for canonical systems},
J. Reine Angew. Math. \textbf{736} (2018), 285--315.

\bibitem{LangerPrucknerWoracek2024}
M.~Langer, R.~Pruckner, and H.~Woracek,
\emph{Estimates for the Weyl coefficient of a two-dimensional canonical system},
Ann. Sc. Norm. Super. Pisa Cl. Sci. (5) \textbf{25} (2024), no.~4, 2259--2330.

\bibitem{LangerWoracek2011}
M.~Langer and H.~Woracek,
\emph{A local inverse spectral theorem for Hamiltonian systems},
Inverse Problems \textbf{27} (2011), no.~5, 055002.

\bibitem{MakarovPoltoratski2023}
N.~Makarov and A.~Poltoratski,
\emph{Etudes for the inverse spectral problem},
J. Lond. Math. Soc. (2) \textbf{108} (2023), no.~3, 916--977.

\bibitem{OrtegaCerdaSeip2002}
J.~Ortega-Cerd\`a and K.~Seip,
\emph{Fourier frames},
Ann. of Math. (2) \textbf{155} (2002), no.~3, 789--806.

\bibitem{PoltoratskiZhang2023}
A.~Poltoratski and A.~R.~Zhang,
\emph{Periodic approximations in inverse spectral problems for canonical Hamiltonian systems},
J. Funct. Anal. \textbf{284} (2023), no.~11, 109883.

\bibitem{Reiffenstein2023}
J.~Reiffenstein,
\emph{A quantitative formula for the imaginary part of a Weyl coefficient},
J. Spectr. Theory \textbf{13} (2023), no.~2, 555--591.

\bibitem{Remling2018}
C.~Remling,
\emph{Spectral Theory of Canonical Systems},
De Gruyter Studies in Mathematics, De Gruyter, Berlin, 2018.

\bibitem{Romanov2014}
R.~Romanov,
\emph{Canonical Systems and de Branges Spaces},
preprint, arXiv:1408.6022, 2014.

\bibitem{Winkler2000}
H.~Winkler,
\emph{Small perturbations of canonical systems},
Integral Equations Operator Theory \textbf{38} (2000), no.~2, 222--250.

\bibitem{Zhang2025}
A.~R.~Zhang,
\emph{Direct spectral problems for Paley--Wiener canonical systems},
preprint, arXiv:2505.00669, 2025.

\end{thebibliography}
\end{document}